\magnification=1200
\overfullrule=0pt
\centerline {\bf On the existence and uniqueness of minima and maxima on spheres}
\centerline {\bf of the integral functional of the calculus of variations}\par
\bigskip
\bigskip
\centerline {BIAGIO RICCERI}\par
\bigskip
\bigskip
{\bf Abstract:} Given a bounded domain $\Omega\subset {\bf R}^n$, we prove that
 if $f:{\bf R}^{n+1}\to {\bf R}$ is a $C^1$ function
whose gradient is Lipschitzian in ${\bf R}^{n+1}$ and non-zero at $0$, then,
for each $r>0$ small enough, the restriction of
the integral functional $u\to \int_{\Omega}f(u(x),\nabla u(x))dx$ to the sphere
$\{u\in H^1(\Omega):\int_{\Omega}(|\nabla u(x)|^2+|u(x)|^2)dx=r\}$ has a unique global minimum
and a unique global maximum.\par
\bigskip
{\bf Key words:} Sobolev space; integral functional; minimum; maximum;
sphere; existence; uniqueness.\par
\bigskip
{\bf 2000 Mathematics Subject Classification.}  49K27.

\bigskip
\bigskip
{\bf Introduction}\par
\bigskip
Here and in the sequel, $\Omega\subset {\bf R}^n$ is a bounded domain,
 and $f:{\bf R}^{n+1}\to {\bf R}$ is a $C^1$ function whose
gradient is non-constant and Lipschitzian (with respect to the
Euclidean metric).\par
\smallskip
We will consider the Sobolev space $H^1(\Omega)$ endowed with the norm
$$\|u\|=\left ( \int_{\Omega}(|\nabla u(x)|^2+|u(x)|^2)dx\right ) ^{1\over 2}\
$$
which is induced by the scalar product
$$\langle u,v\rangle=\int_{\Omega}(\nabla u(x)\nabla v(x)+u(x)v(x))dx\ .$$
\smallskip
The linear growth of $\nabla f$ (coming from its Lipschitzianity) implies
 that
the functional
$$u\to J(u):=\int_{\Omega}f(u(x),\nabla u(x))dx$$
is (well-defined and) $C^1$ on $H^1(\Omega)$, with derivative given by
$$\langle J'(u),v\rangle=\int_{\Omega}(f_{\xi}(u(x),\nabla u(x))v(x)+
\nabla_{\eta} f(u(x),\nabla u(x))\nabla v(x))dx$$
for all $u,v\in H^1(\Omega)$ ([2], p. 249).\par
\smallskip
Let $r>0$. We are interested in minima and maxima of the restriction
of the functional $J$ to
the sphere $S_r:=\{u\in H^{1}(\Omega) :\|u\|=r\}$.\par
\smallskip
In the present setting, there is no evidence of their existence and
uniqueness. In
fact, with regard to the existence aspect, not only $S_r$ is not weakly
compact but also, if $f(\xi,\cdot)$ is neither
convex nor concave in ${\bf R}^n$, the functional $J$ is neither lower nor upper weakly
semicontinuous. But, even when $J$ is sequentially weakly continuous, it
may happen that $J$ has no minima and/or maxima on $S_r$.
\smallskip
In this connection, consider the following simple and enlightenting
situation. Assume that $f$ depends only on the first variable and that
has a unique global
maximum in ${\bf R}$, say $\xi_0$. So, $J(u)=\int_{\Omega}f(u(x))dx$.
Then, it is clear that the constant function $x\to \xi_0$ is the unique
maximum of the functional $J$. In this case, $J$ turns out to be sequentially
weakly continuous, thanks to the Rellich-Kondrachov theorem. Then,
by Lemma 2.1 of [1], the function $\rho\to \sup_{S_{\rho}}J$ is non-decreasing
in $]0,+\infty[$. Consequently, if $r>|\xi_0|(\hbox {\rm meas}(\Omega))^
{1\over 2}$, $J_{|S_r}$ has no maxima.\par
\smallskip
Nevertheless, we will show that if $\nabla f(0)\neq 0$ then
$J_{|S_r}$ possesses exactly one minimum and exactly one maximum for
each $r>0$ small enough.\par
\bigskip
{\bf The result} \par
\bigskip
To shorten the statement of our result, let us introduce some further notations.
In the sequel, $g:{\bf R}^{n+1}\to {\bf R}$ is another $C^1$ function which
is non-negative, with $g(0)=0$, and whose gradient
is Lipschitzian, with Lipschitz constant $\nu<2$. We set
$$I(u)=\int_{\Omega}g(u(x),\nabla u(x))dx$$
for all $u\in H^1(\Omega)$.\par
\smallskip
Moreover, $V$ is a closed linear subspace of $H^1(\Omega)$ with the following
property: there exists $v_0\in V$ such that
$$\int_{\Omega}(f_{\xi}(0)v_0(x)+\nabla_{\eta}f(0)\nabla v_0(x))dx\neq 0\ .$$
\indent
Finally, if $L$ is the Lipschitz constant of
$\nabla f$, we denote by $S$ the set (possibly empty) of all global minima
of the restriction to $V$ of the functional
$$u\to \|u\|^2+I(u)+{{2-\nu}\over {L}}J(u)\ .$$
\indent
Then, with the convention $\inf\emptyset=+\infty$, our result reads as follows:\par
\vfill\eject
THEOREM 1. - {\it Under the above assumptions, one has\par
$$\delta:=\inf_{u\in S}(\|u\|^2+I(u))>0$$
 and, for each $r\in ]0,\delta[$, the restriction
of the functional $J$ to the set $$C_r:=\left \{ u\in V :\|u\|^2+I(u)=r\right \}$$
has a unique global minimum.}\par
\smallskip
PROOF. Let $\mu\geq 0$ and let $u, v, w\in H^{1}(\Omega)$, with $\|w\|=1$.
Using Cauchy-Schwartz and H\"older inequalities, we have
$$|\langle I'(u)+\mu J'(u))-I'(v)-\mu J'(v),w\rangle|\leq$$
$$\leq \int_{\Omega}
|(g_{\xi}(u,\nabla u)-g_{\xi}(v,\nabla v))w+
(\nabla_{\eta} g(u,\nabla u)-\nabla_{\eta} g(v,\nabla v))\nabla
w|dx+$$
$$+\mu
\int_{\Omega}|(f_{\xi}(u,\nabla u)-f_{\xi}(v,\nabla v))w+
(\nabla_{\eta} f(u,\nabla u)-\nabla_{\eta} f(v,\nabla v))\nabla
w|dx\leq$$
$$\leq
\int_{\Omega}\left (
|(g_{\xi}(u,\nabla u)-g_{\xi}(v,\nabla v)|^2+
|\nabla_{\eta} g(u,\nabla u)-\nabla_{\eta} g(v,\nabla v)|^2\right ) ^{1\over 2}
\left ( |w|^2+|\nabla w|^2\right ) ^{1\over 2}dx+$$
$$+\mu
\int_{\Omega}\left ( 
|(f_{\xi}(u,\nabla u)-f_{\xi}(v,\nabla v)|^2+
|\nabla_{\eta} f(u,\nabla u)-
\nabla_{\eta} f(v,\nabla v)|^2\right ) ^{1\over 2}
\left ( |w|^2+|\nabla w|^2\right ) ^{1\over 2}dx\leq$$
$$\leq \left ( \int _{\Omega}
(|(g_{\xi}(u,\nabla u)-g_{\xi}(v,\nabla v)|^2+
|\nabla_{\eta} g(u,\nabla u)-\nabla_{\eta} g(v,\nabla v)|^2)dx\right )
 ^{1\over 2}+$$
$$+\mu \left ( \int _{\Omega}
(|(f_{\xi}(u,\nabla u)-f_{\xi}(v,\nabla v)|^2+
|\nabla_{\eta} f(u,\nabla u)-\nabla_{\eta} f(v,\nabla v)|^2)dx\right )
 ^{1\over 2}\leq$$
$$\leq (\nu+\mu L)\|u-v\|\ .$$
Hence, the derivative of the functional $I+\mu J$ is
Lipschitizian, with constant $\nu+\mu L$. As a consequence,
if $0\leq \mu <{{2-\nu}\over {L}}$, the functional
$u\to \|u\|^2+I(u)+\mu J(u)$ is strictly convex and coercive. To see this, it is enough
to show that its derivative is strongly monotone ([3], pp. 247-248).
Indeed, if $\Phi(\cdot):=\|\cdot\|^2$, 
we have for all $u, v\in H^1(\Omega)$
$$\langle\Phi'(u)+I'(u)+\mu
J'(u)-\Phi'(v)-I'(v)-\mu J'(v),u-v\rangle\geq$$
$$\geq 2\|u-v\|^2-\|I'(u)-J'(v)+\mu (I'(v)-J'(v))\|\|u-v\|\geq
(2-\nu-\mu L)\|u-v\|^2\ .$$
Clearly, this shows also the convexity of the functional
$\Phi+I+{{2-\nu}\over {L}}J$. Assume $S\neq \emptyset$. Then, $S$ is
closed and convex, and so there exists a unique $\hat u\in S$ such that
$$\|\hat u\|^2+I(\hat u)=\delta\ .$$
Observe that $\|u\|^2+I(u)>0$ for all $u\in V\setminus \{0\}$. So, $\delta\geq
0$. Arguing by contradiction, assume $\delta=0$. Then, it would
follow $\hat u=0$. Hence, since $0\in S$, we would have
$$\langle\Phi'(0)+I'(0)+{{2-\nu}\over {L}}J'(0),v\rangle=0$$
for all $u\in V$
and so, since $\Phi'(0)+I'(0)=0$ (being $0$ the global minimum of $\Phi+I$),
it would follow
$$\int_{\Omega}(f_{\xi}(0)v(x)+
\nabla_{\eta} f(0)\nabla v(x))dx=0$$
for all $v\in V$, against one of the hypotheses. Hence, we have proven that $\delta>0$.
Now, fix $r\in ]0,\delta[$ and consider the function $\Psi:V\times
[{{L}\over {2-\nu}},+\infty[\to {\bf R}$ defined by
$$\Psi(u,\lambda)=J(u)+\lambda (\|u\|^2+I(u)-r)$$
for all $(u,\lambda)\in V\times [{{L}\over {2-\nu}},+\infty[)$.
As we have seen above, $\Psi(\cdot,\lambda)$ is continuous and convex for
all $\lambda\geq {{L}\over {2-\nu}}$ and coercive for all $\lambda>
{{L}\over {2-\nu}}$, while $\Psi(u,\cdot)$ is continous and concave for
all $u\in V$, with $\lim_{\lambda\to +\infty}\Psi(0,\lambda)=-\infty$.
So, we can apply to $\Psi$ a classical saddle-point theorem
([3], Theorem 49.A) which
ensures the existence of $(u^*,\lambda^*)\in V\times [{{L}\over {2-\nu}},+\infty[$
such that
$$J(u^*)+\lambda^* (\|u^*\|^2+I(u^*)-r)=\inf_{u\in V}(J(u)+\lambda^* (\|u\|^2+I(u)-r))=$$
$$=J(u^*)+\sup_{\lambda\geq {{L}\over {2-\nu}}}\lambda (\|u^*\|^2+I(u^*)-r)\ .$$
Of course, we have $\|u^*\|^2+I(u^*)\leq r$, since the sup is finite. But, if it were
$\|u^*\|^2+I(u^*)<r$, we would have $\lambda^*={{L}\over {2-\nu}}$. This, in turn, would
imply that $u^*\in S$, against the fact that $r<\delta$. Hence, we have
$\|u^*\|^2+I(u^*)=r$. Consequently
$$J(u^*)+\lambda^* r=\inf_{u\in V}(J(u)+\lambda^* (\|u\|^2+I(u)))\ .$$
 From this, we infer that $\lambda^*>{{L}\over {2-\nu}}$ (since $r<\delta$), that
$u^*$ is a global minimum of $J_{|C_r}$ and that if each global minimum of
$J_{|C_r}$ is a global minimum in $V$ of the functional $u\to \|u\|^2+I(u)+
\lambda^* J(u)$. Since $\lambda^*>{{L}\over {2-\nu}}$, this functional is strictly
convex and so $u^*$ is its unique global minimum in $V$. The proof is complete.\hfill
$\bigtriangleup$
\medskip
REMARK 1. It is almost superfluous to remark that the conclusion
of Theorem 1 may fail if the assumption that involves $V$ and $\nabla f(0)$ is
not satisfied. In this connection, consider, for instance,
the case $f(\sigma)=-|\sigma|^2$, with $g=0$. This assumption, however,
serves only to ensure that $\delta>0$. So, it becomes superfluous, in particular,
when $S=\emptyset$.
\par
\medskip
Now, denote by $S_1$ the set (possibly empty) of all global minima
of the restriction to $V$ of the functional
$$u\to \|u\|^2+I(u)-{{2-\nu}\over {L}}J(u)\ .$$
Clearly, applying Theorem 1 also to $-f$, we get\par
\medskip
THEOREM 2. - {\it Under the assumptions of Theorem 1, one has
$$\delta_{1}:=\min\left \{ \inf_{u\in S}(\|u\|^2+I(u)),\inf_{u\in S_1}(\|u\|^2+I(u))\right \} 
>0$$
 and, for each $r\in ]0,\delta_{1}[$, the restriction
of the functional $J$ to the set $$\left \{ u\in V :\|u\|^2+I(u)=r\right \}$$
has a unique global minimum and a unique global maximum.}\par
\bigskip
\bigskip
\bigskip
\bigskip
\centerline {\bf References}\par
\bigskip
\bigskip
\noindent
[1]\hskip 5pt M. SCHECHTER and K. TINTAREV, {\it Spherical maxima in
Hilbert space and semilinear elliptic eigenvalue problems}, Differential
Integral Equations, {\bf 3} (1990), 889-899.\par
\smallskip
\noindent
[2]\hskip 5pt M. STRUWE, {\it Variational methods}, Springer-Verlag, 1996.\par
\smallskip
\noindent
[3]\hskip 5pt  E. ZEIDLER, {\it Nonlinear functional analysis and
its applications}, vol. III, Springer-Verlag, 1985.\par 
\bigskip
\bigskip
\bigskip
\bigskip
Department of Mathematics\par
University of Catania\par
Viale A. Doria 6\par
95125 Catania\par
Italy\par
{\it e-mail address}: ricceri@dmi.unict.it

\bye